\documentclass[12pt]{amsart}
\usepackage{amsfonts}

\usepackage{amssymb}
\usepackage{enumitem}

\makeatletter
\@namedef{subjclassname@2010}{  \textup{2010} Mathematics Subject Classification}
\makeatother

\theoremstyle{definition}

\numberwithin{equation}{section}
\numberwithin{equation}{section}

\begin{document}
\title[On Sj\"{o}lin-Soria-Antonov type extrapolation for  locally compact groups]{On Sj\"{o}lin-Soria-Antonov type extrapolation for  locally compact groups}
\author[G. Oniani]{ Giorgi Oniani}

\address{Department of Mathematics \\
Akaki Tsereteli State University \\
59 Tamar Mepe St., Kutaisi 4600\\
Georgia}
\email{oniani@atsu.edu.ge}
\date{}

\begin{abstract}
Sj\"{o}lin-Soria-Antonov type extrapolation theorem for  locally compact $\sigma$-compact non-discrete  groups is proved. As an application of this result it is shown that the Fourier series with respect to the Vilenkin orthonormal systems on the Vilenkin groups of bounded type   converge almost everywhere for functions from the class $L\log^{+}L\log^{+}\log^{+}\log^{+}L$.

\end{abstract}

\maketitle

\baselineskip=17pt

Let $(X,\mu)$ be a measure space. Denote by:
\begin{enumerate}
\item[$\bullet$] $L^0(X,\mu)$ the class of all measurable functions $f:X\rightarrow [-\infty,\infty]$;

\item[$\bullet$]  $\Phi$ the set of all increasing continuous functions $\varphi:[0,\infty)\rightarrow[0,\infty)$ with $\varphi(0)=0$ and  $\liminf\limits_{u\rightarrow\infty}\varphi(u)/u>0$;

\item[$\bullet$] $\varphi(L)(X,\mu)$  the class of all measurable functions $f:X\rightarrow [-\infty,\infty]$ for which $\int_X \varphi(|f|)d\mu<\infty$;

\item[$\bullet$]  $\chi_E$ the characteristic function of a set $E\subset X$.

\end{enumerate}

We say that an operator $T: L(X,\mu)\rightarrow L^0(X,\mu)$ is:
\begin{enumerate}
\item[$\bullet$] \emph{of restricted weak type $\varphi$}, where $\varphi \in\Phi$, if there is a number  $C>0$ such that $\mu(\{|T \chi_E|>\lambda\})\leq C\varphi (1/\lambda) \mu(E)$
    for every $E\subset X$ with finite measure and every $\lambda>0$;

\item[$\bullet$] \emph{of restricted weak type $p$}, where $1\leq p< \infty$, if $T$ is of weak type $\varphi$ for the case of the function $\varphi(u)=u^p$;

\item[$\bullet$] \emph{of generalized weak type $\varphi$}, where $\varphi \in\Phi$, if for every subset $Y\subset X$  with finite measure and every $\varepsilon>0$ there is a number $C_{Y,\varepsilon}>0$ such that $\mu(\{|T f|>\lambda\}\cap Y)\leq \varepsilon +C_{Y,\varepsilon}\int_X\varphi (|f|/\lambda)d\mu$ for every $f\in L(X,\mu)$ and every $\lambda>0$.

\end{enumerate}

Denote by:
\begin{enumerate}
\item[$\bullet$] $\mathrm{H}_m(X,\mu)$ $(m>0)$ the class of all sub-linear  operators $T:L(X,\mu) \rightarrow L^0(X,\mu)$ for which  there are $r=r(T)>1$ and $C=C(T)>0$ such that $T$ is of restricted weak  type $p$ for every $p\in (1,r)$ with constant $C_p$ satisfying the estimation  $C_p\leq \big(\frac{C}{p-1} \big)^{mp}$;

\item[$\bullet$]  $\mathrm{A}_m(X,\mu)$ $(m>0)$ the class of all sub-linear operators $T:L(X,\mu) \rightarrow L^0(X,\mu)$ which are of generalized weak type $\varphi$, where $\varphi(u)=u(1+\log^{+}u)^m$ $(1+\log^{+}\log^{+}\log^{+}u)$;

\end{enumerate}

For a family $\mathbf{T}=(T_j)$  of operators acting from  $L(X,\mu)$ to  $L^0(X,\mu)$  \emph{ the maximal operator} $M_{\mathbf{T}}$ is defined as follows:
$$
M_{\mathbf{T}} f(x)=\sup_{j} |T_j f (x)|\;\;\;\;(f\in L(X,\mu), x\in X).
$$

\textbf{Remark 1}.  If a sequence $(T_j f)$ converge a.e. for each function $f$ from some "dense" subclass of $\varphi(L)$, then the generalized weak type $\varphi$ estimation for $M_{\mathbf{T}}$ by standard technique (see, e.g., Lemma 2 below) makes it possible to establish the a.e. convergence for each function $f$ from $\varphi(L)$.

Let $X$  be  a locally compact group and $\mu$   be a left-invariant  Haar measure in $X$. For a function $k\in L(X,\mu)$ denote by $T_k$ \emph{the convolution operator generated by} $k$, i.e., $T_k f=f\ast k$ $(f\in L(X,\mu))$. Denote by $\mathrm{M}(X,\mu)$ the class of all  maximal operators $M_{\mathbf{T}}$, where $\mathbf{T}$ is a sequence of convolution operators $(T_{k_j})$ for some functions   $k_j\in L(X,\mu)$ $(j\in\mathbb{N})$.

 \emph{ A  Vilenkin group } is defined as  the  direct product of discrete cyclic groups $\mathbb{Z}_{m_j}=\{0,1,\dots, m_j-1\}$   with $m_j\geq 2$ $(j\in\mathbb{N})$. A Vilenkin group is said to be of \emph{bounded type} if $\sup_j m_j<\infty$.  By $(\xi^X_j)$ we will denote \emph{the Vilenkin orthonormal system } corresponding to a Vilenkin group $X$ (see, e.g., [1, App. (0.7)] or  [2, \S 1.5] for the definition). Note that for the case $m_j=2$ $(j\in\mathbb{N})$ the terms \emph{dyadic group} and \emph{Walsh orthonormal system} respectively are used.

Below $\mathbb{T}$ denotes the one-dimensional torus, i.e., $\mathbb{T}=\mathbb{R}/\mathbb{Z}$. Throughout the paper we will use the following convention: $\log n$  stands for $\log_2 n$.


\section{Results}

According to the famous theorem  of Carleson [3] the Fourier series of every function $f\in L^2(\mathbb{T})$ converge almost everywhere.  Hunt [4] extended this result to the spaces $L^p(\mathbb{T})$ $(p>1)$, furthermore, in [4] it was proved that the  Carleson maximal operator (i.e. the operator $M_{\mathbf{T}}$, where $\mathbf{T}=(T_j)$  is the sequence of Fourier partial sums) belongs to the class $\mathrm{H}_1(X,\mu)$. Using this  estimation  Antonov [5] have shown  almost everywhere convergence of  Fourier series of functions from the class  $L\log^{+}L\log^{+}\log^{+}\log^{+}L(\mathbb{T})$. Note that on the other hand, due to the result of Konyagin [6] in classes $\varphi(L)$ with $\varphi(u)=o(u\sqrt{\log u/ \log\log u})$ $(u\rightarrow\infty)$  it is not guaranteed a.e. convergence of the Fourier series.
Refining the method used in [5] Sj\"{o}lin and Soria [7] proved the following extrapolation principle: Let  $(X,\mu)$ be  either  $\mathbb{R}^n$ or $\mathbb{T}^n$ with Lebesgue measure and $m>0$. Then $\mathrm{M}(X,\mu)\cap \mathrm{H}_m(X,\mu)\subset \mathrm{A}_m(X,\mu)$.  The following extension of this result is true.

\textbf{Theorem 1}. \emph{Let $X$  be  a locally compact $\sigma$-compact non-discrete group,  $\mu$   be a left-invariant  Haar measure in $X$ and $m>0$. Then}
$$
\mathrm{M}(X,\mu)\cap \mathrm{H}_m(X,\mu)\subset \mathrm{A}_m(X,\mu). \eqno(1)
$$

\textbf{Remark 2}.  Theorem 1 implies the validity of the inclusion $(1)$ for the cases when $X$  is either locally compact non-discrete group satisfying the second axiom of countability or compact non-discrete group.

Applying Theorem 1 to the a.e. convergence problem of Fourier series with respect to orthonormal systems of characters of compact groups we prove the following theorems.

\textbf{Theorem 2}. \emph{Let $X$ be a compact Abelian group, $\mu$ be a  Haar measure in $X$, $(\xi_j)$ be a complete orthonormal system consisting of characters of $X$, $\mathbf{T}=(T_j)$ be  the sequence of partial sums of a Fourier series with respect to $(\xi_j)$ and $m>0$.  If the maximal operator $M_{\mathbf{T}}$ belongs to the class $\mathrm{H}_m(X,\mu)$ then the Fourier series with respect  to $(\xi_j)$ of every  function $f$ from the class  $L(\log^{+}L)^m\log^{+}\log^{+}\log^{+}L(X,\mu)$ converge to $f$ almost everywhere.}

\textbf{Theorem 3}. \emph{Let $X$ be a  Vilekin group of bounded type  and $\mu$ be a Haar measure in $X$. Then the Fourier series with respect  to the Vilenkin system $(\xi^X_j)$ of every  function $f$ from the class  $L\log^{+}L\log^{+}\log^{+}\log^{+}L(X,\mu)$ converge to $f$ almost everywhere.}

In the the setting of Theorem 3 a.e. convergence for functions from the classes $L^p(X,\mu)$ $(p>1)$ was proved by Gosselin [8].

In [7] it was shown a.e. convergence of Walsh-Fourier series for functions from $L\log^{+}L\log^{+}\log^{+}\log^{+}L([0,1])$. This result is equivalent to the conclusion of Theorem 3  for the case of the Walsh orthonormal system  on the dyadic group.
It must be mentioned here that by virtue of the  result of Bochkarev [9] in classes $\varphi(L)$ with  $\varphi(u)=o(u\sqrt{\log u})$ $(u\rightarrow\infty)$ it is not guaranteed a.e. convergence of the Walsh-Fourier series.  Analogues theorem for the  Vilenkin systems on Vilenkin groups of bounded type was established by Polyakov [10].


\section{Proofs}

The Lemma 1 below extends the approximation principle proved in [7] (see Lemma 5) to the setting of Theorem 1. The rest part of the proof of Theorem 1 follows the scheme used in [7].

\textbf{Lemma 1}.  \emph{Let $X$ be a locally compact $\sigma$-compact non-discrete  group and $\mu$ be a left-invariant Haar measure in $X$. Suppose $(k_j)_{j\geq1}$ is a sequence of functions from $L(X,\mu)$. Then for every $N\in\mathbb{N}$, a non-negative function $f\in L(X,\mu)$, increasing sequence $(a_n)_{n\geq 0}$ with $a_0=0$ and $\lim\limits_{n\rightarrow\infty} a_n=\infty$, and   $\varepsilon>0$,  there is a simple function $h$ of the form $h=\sum_{n=1}^{\nu}a_n \chi_{E_n}$ for some $\nu$ such that}
\begin{enumerate}

\item[$\bullet$] $E_n\subset\{a_{n-1}<f\leq a_n\}$ \emph{and }$\int_{\{a_{n-1}<f\leq a_n\}}f d\mu=a_n \mu(E_n)$ \emph{for every} $n=1,\dots, \nu$;

\item[$\bullet$] $\int_{X} M_{(T_{k_1},\dots, T_{k_N})}(f-h)d\mu < \varepsilon$.

\end{enumerate}

\textbf{Proof}. For the sake of simplicity we will consider the case when $X$ is Abelian.

Denote   $f_t=f\chi_{\{f<t\}}$ and  $f^t=f\chi_{\{f\geq t\}}$ $(t>0)$. For every $t>0$ we have $$
\int_X M_{(T_{k_1},\dots, T_{k_N})}(f^t)d\mu\leq\sum_{j=1}^N \int_X |f^t\ast k_j|d\mu\leq \sum_{j=1}^N \|f^t\|_{L(X,\mu)} \| k_j\|_{L(X,\mu)}.
$$
Since $\lim\limits_{t\rightarrow\infty}\|f^t\|_{L(X,\mu)}=0$ then by the above estimation we can choose $\nu$ so big that
$$
\int_X M_{(T_{k_1},\dots, T_{k_N})} (f^{a_\nu})d\mu<\frac{\varepsilon}{2}. \eqno(2)
$$

 For every function $f\in L(X,\mu)$ the mapping $X\ni s\mapsto f(\cdot - s)\in L(X,\mu)$ is uniformly continuous (see, e.g., [11, (20.4) Theorem]). Using this statement for functions $k_1,\dots, k_N$ we find a neighbourhood $V$ of zero in $X$  such that
 $$
 \int_X |k_j(x-s_2)-k_j(x-s_1)|d\mu(x)<\frac{\varepsilon}{4N(1+\|f\|_{L(X,\mu)})} \eqno (3)
 $$
for every $j=1,\dots, N$ and $s_1, s_2\in X$ with $(s_2-s_1)\in V$.

Let $X_m$ $(m\in\mathbb{N})$ be compact subsets of $X$ such that $X=\bigcup_{m=1}^\infty X_m$.  For each $m$ we can find finite set $X_m^\ast \subset X_m$ for which $X_m\subset\bigcup_{x\in X_m^\ast}(x+V)$. Then the collection $\{x+V: x\in\bigcup_{m=1}^\infty X_m^\ast \}$ is at most countable. Without loss of generality let us assume that this collection is countable and let $V_1, V_2, \dots $ be its members. Set $V_1^\ast=V_1$ and  $V_m^\ast=V_m\setminus \bigcup_{i=1}^{m-1}V_i$ when $m\geq 2$. Note that $V_m^\ast\subset V_m$  $(m\in\mathbb{N})$, $V_{m_1}^\ast\cap V_{m_2}^\ast=\emptyset$ $(m_1\neq m_2)$ and  $X=\bigcup_{m=1}^\infty V_m^\ast$.  Without loss of generality assume also that $\mu(V_m^\ast)>0$ for each $m$.

Denote $G_n=\{a_{n-1}<f\leq a_n\}$ $(n\in\mathbb{N})$. For every $n=1,\dots,\nu$ and $m\in\mathbb{N}$ we have
$$
\int_{G_n\cap V_{m}^\ast}fd\mu\leq a_n\mu(G_n\cap V_{m}^\ast).
$$
Since the group $X$ is non-discrete then the measure $\mu$ is non-atomic (see, e.g., [12, \S 58]). Consequently, by virtue of Liapounoff's theorem (see, e.g., [13, Theorem 5.5]) there is a set $E_{n,m}\subset G_n\cap V_{m}^\ast$ such that
$$
\int_{G_n\cap V_{m}^\ast}fd\mu=a_n \mu(E_{n,m}).\eqno(4)
$$
Denote $E_n=\bigcup_{m=1}^\infty E_{n,m}$ $(n=1,\dots, \nu)$ and $h=\sum_{n=1}^{\nu}a_n\chi_{E_n}$. Then
$$
\int_X M_{(T_{k_1},\dots, T_{k_N})}(f_{a_\nu}-h)d\mu\leq \sum_{j=1}^N \int_X |(f-h)\ast k_j|d\mu\leq
$$
$$
 \sum_{j=1}^N \sum_{n=1}^\nu \sum_{m=1}^\infty \int_X \bigg|\int_{G_n\cap V_{m}^\ast} k_j(x-y)(f(y)-a_n\chi_{E_{n,m}}(y))d\mu(y)\bigg|d\mu(x).\eqno(5)
$$
By virtue of $(4)$, for every $n=1,\dots,\nu$ and $m\in\mathbb{N}$ we have
$$
\int_{G_n\cap V_{m}^\ast}(f(y)- a_n\chi_{E_{n,m}}(y))d\mu(y)=0. \eqno(6)
$$
Let for every $m\in\mathbb{N}$, $x_m$ be the  point for which $V_m=x_m+V$. Now using  $(3-6)$ we obtain that
$$
\int_X M_{(T_{k_1},\dots, T_{k_N})}(f_{a_\nu}-h)d\mu\leq
$$
$$
 \sum_{j,n,m} \int\limits_X \bigg|\int\limits_{G_n\cap V_{m}^\ast} (k_j(x-y)-k_j(x-x_m))(f(y)-a_n\chi_{E_{n,m}}(y))d\mu(y)\bigg|d\mu(x)\leq
$$
$$
 \sum_{j,n,m} \; \int\limits_{G_n\cap V_{m}^\ast} (f(y)+a_n\chi_{E_{n,m}}(y))\bigg(\int\limits_X |k_j(x-y)-k_j(x-x_m)|d\mu(x)\bigg)d\mu(y)
$$
$$
\leq \sum_{j=1}^N \sum_{n=1}^\nu \sum_{m=1}^\infty\frac{\varepsilon}{4N(1+\|f\|_{L(X,\mu)})}\int\limits_{G_n\cap V_{m}^\ast} (f(y)+a_n\chi_{E_{n,m}}(y))d\mu(y)=
$$
$$
\frac{\varepsilon}{4(1+\|f\|_{L(X,\mu)})}\sum_{j=1}^N 2\int\limits_{G_n} f(y)d\mu(y)\leq\frac{\varepsilon \|f\|_{L(X,\mu)}}{2(1+\|f\|_{L(X,\mu)})} <\frac{\varepsilon}{2}. \eqno(7)
$$
Combining the estimations $(2)$ and $(7)$ and using sub-linearity of the operator $M_{(T_{k_1},\dots, T_{k_N})}$ we conclude the validity of the lemma.

\medskip

Let $(X,\mu)$ be a measure space with finite measure and $\varphi\in\Phi$.  We say that a class of measurable functions  $\Delta$ is  $\varphi$-\emph{dense in} $\varphi (L) (X,\mu)$ if for every $f\in \varphi (L) (X,\mu)$ and $\varepsilon>0$ there is a function $h\in \Delta$ such that $\int_X \varphi(|f - h|)d\mu<\varepsilon.$ Note that using the dominated convergence theorem it is easy to check $\varphi$-density in $\varphi (L) (X,\mu)$ of the class of all bounded measurable functions.

A function $\varphi:[0,\infty)\rightarrow [0,\infty)$ is said to satisfy $\Delta_2$-\emph{condition} if there are positive numbers $c$ and $u_0$ such that $\varphi(2u)\leq c\varphi(u)$ when $u\geq u_0$.

\textbf{Lemma 2}. \emph{Let $(X,\mu)$ be a measure space with finite measure, $\mathbf{T}=(T_j)$ be a sequence of linear operators acting from $L(X,\mu)$ to $L^0(X,\mu)$,   $\Pi$ be the class of all functions $f\in L(X,\mu)$ for which $(T_j f)$ converge to $f$ almost everywhere,  and $\varphi \in\Phi$ be a function satisfying the $\Delta_2$-condition. If $\;\Pi$ is $\varphi$-dense in $\varphi (L) (X,\mu)$ and $M_{\mathbf{T}}$ is of generalized weak type $\varphi$ then $\varphi (L) (X,\mu)\subset \Pi$.}

 \textbf{Proof}. For   $f\in L(X,\mu)$ and $\lambda>0$ denote
 $$
 E(f,\lambda)=\{\limsup\limits_{j\rightarrow\infty}|T_j f-f|>\lambda\}.
 $$
Note that if $f\in \Pi$ then  $\mu(E(f,\lambda))=0$ for every $\lambda>0$.

Let $f,h\in L(X,\mu)$. Since   for every $x\in X$,
$$
\limsup\limits_{j\rightarrow\infty}|T_j f(x)-f(x)|\leq
$$
$$
\leq\limsup\limits_{j\rightarrow\infty}|T_j f(x)-T_j h(x)|+ \limsup\limits_{j\rightarrow\infty}|T_j h(x)-h(x)| + |f(x)- h(x)|
$$
then  for every $\lambda>0$ we have
$$
E(f,\lambda)\subset \{\limsup\limits_{j\rightarrow\infty}|T_j f-T_j h|>\lambda/3\}\cup E(h,\lambda/3)\cup \{|f-h|>\lambda/3\}\subset
$$
$$
\subset \{M_{\mathbf{T}}(f-h)>\lambda/3\}\cup E(h,\lambda/3)\cup \{|f-h|>\lambda/3\}. \eqno(8)
$$

 Suppose, $f\in \varphi(L)(X,\mu)$. Let us show that for every  $\lambda\in (0,1)$ the set $E(f,\lambda)$ has zero measure. It will imply the needed conclusion.  Let $\varepsilon$ be an arbitrary positive number. For every $h\in L(X,\mu)$ we have
 $$
 \mu(\{M_{\mathbf{T}}(f-h)>\lambda/3\})\leq \varepsilon+ C_{X,\varepsilon}\int_X \varphi\bigg(\frac{3|f-h|}{\lambda}\bigg)d\mu, \eqno(9)
 $$
 and
 $$
 \mu(\{|f-h|>\lambda/3\})\leq \frac{1}{\varphi(\lambda/3)}\int_X \varphi(|f-h|)d\mu.\eqno(10)
 $$

 Since $\varphi$ satisfy the $\Delta_2$-condition then $3f/\lambda \in \varphi(L)(X,\mu)$. Consequently, we can choose a function $g\in \Pi$ so that
 $$
 \int_X\varphi\bigg(\bigg|\frac{3f}{\lambda}-g\bigg|\bigg)d\mu<\frac{\varepsilon}
 {C_{X,\varepsilon}+1/\varphi(\lambda/3)}. \eqno(11)
 $$

Set $h=\lambda g/3$. Note that the class $\Pi$ is  positively homogeneous. Therefore $h\in\Pi$ and consequently,
$$
\mu(E(h,\lambda/3))=0.\eqno(12)
$$
On the other hand we have
$$
\int_X\varphi\bigg(\frac{3|f-h|}{\lambda}\bigg)d\mu=\int_X\varphi\bigg(\bigg|\frac{3f}{\lambda}-
g\bigg|\bigg)d\mu \eqno(13)
$$
and
$$
\int_X \varphi(|f-h|)d\mu=\int_X\varphi\bigg(\frac{\lambda}{3}\bigg|\frac{3f}{\lambda}-
g\bigg|\bigg)d\mu\leq \int_X\varphi\bigg(\bigg|\frac{3f}{\lambda}-
g\bigg|\bigg)d\mu. \eqno(14)
$$
 From $(8-14)$ we obtain the estimation $\mu(E(f,\lambda))<3\varepsilon$.  Consequently, by arbitrariness of $\varepsilon>0$ we have $\mu(E(f,\lambda))=0$. The lemma is proved.

\textbf{Remark 3}.  Lemma 2 remains true without assumption of $\varphi$ to satisfy the $\Delta_2$-condition. In this setting the proof becomes a bit longer because there appears an additional step. Namely, first we have to approximate a function $f\in\varphi(L)(X,\mu)$  by its lower truncation $f_t=f\chi_{\{f<t\}}$ and then $f_t$ by a function $h\in \Pi$.

\textbf{Proof of Theorem 2}.  Note that the  partial sum operators $T_j$ of a Fourier series with respect to the system of characters $(\xi_j)$ is represented by the convolution (in $X$) with Dirichlet kernels corresponding to $(\xi_j)$, i.e. $T_j(f)=f\ast D_j$, where $D_j=\sum_{i=0}^{j-1} \xi_i$. Consequently, by virtue of Theorem 1 the maximal operator $M_{\mathbf{T}}= M_{(T_{j})}$ is of generalized weak type $\varphi$, where $\varphi(u)=(1+\log^{+}u)^m (1+\log^{+}\log^{+}\log^{+}u)$.

Due to the interpolation theorem of Stein and Weiss (see, e.g., [14, Corollary 1.4.21]), $M_{\mathbf{T}}$ is of strong type $p$ for every $p\in (1, r(M_{\mathbf{T}}))$. Consequently, assuming  $p\in (1, \min\{r(M_{\mathbf{T}}),2\})$, and, taking into account density in $L^2(X,\mu)$ of the class of all polynomials  with respect to the system  $(\xi_j)$, we conclude that the Fourier series with respect  to $(\xi_j)$ of every  function $f$ from  $L^p(X,\mu)$ converge almost everywhere. Now using   $\varphi$-density in $\varphi(L)(X,\mu)$ of the  class of bounded measurable functions, on the basis  of Lemma 2 we obtain the needed conclusion.

\textbf{Proof of Theorem 3}. Let $\mathbf{T}=(T_j)$ be the sequence of partial sum operators of a Fourier series with respect to the system $(\xi^X_j)$.  Gosselin [8] proved that for every $q\in (1,2)$ there is a constant $C_q>0$ for which
$$
\mu(\{M_{\mathbf{T}} \chi_E>\lambda\})\leq  \frac{C_q}{\lambda^q}\mu(E),
$$
for every measurable set $E\subset X$ and $\lambda>0$. Perusing the proof of this result one can observe that (see the estimation $(84)$ for $\mu(E_1)$, the estimation for $\mu(E_2)$ after the relation $(93)$ and the estimation  for $|S_nf(x)|$ at the end of the proof of the basic result in [8]) for the constants $C_q$  the estimation $C_q\leq \big(\frac{C}{q-1} \big)^{q}$ $(q\in(1,2))$ is valid, where   $C>0$ is an absolute constant. After it by virtue of Theorem 2 we finish the proof.

\end{document}